\documentclass[a4paper,12pt]{amsart}
\usepackage{amscd}
\usepackage{amsmath}
\usepackage{amssymb}

\theoremstyle{plain}
\newtheorem{theorem}{Theorem}[section]
\newtheorem{lemma}[theorem]{Lemma}

\theoremstyle{definition}

\newtheorem{remark}[theorem]{Remark}

\newtheoremstyle{notation}{14pt}{14pt}
     {}
     {}
     {\bfseries}
     {.}
     {\newline}
     {}

\theoremstyle{notation}
\newtheorem*{notation}{Notations}

\input xy
\usepackage[all]{xy}

\begin{document}

\def\pr{\mathbb{P}}
\def\com{\mathbb{C}}
\def\zet{\mathbb{Z}}
\def\real{\mathbb{R}}
\def\field{\mathbb{K}}
\def\rational{\mathbb{Q}}
\def\nat{\mathbb{N}}
\def\derxcat{{\bf D}(X)}
\def\dermcat{{\bf D}(M)}
\def\derccat{{\bf D}(C)}
\def\rhom{R{\mathrm{Hom}}}
\def\gistr{\mathbb{G}}
\def\qiso{\simeq_{\mathrm{q.iso}}}
\def\pic{\mathrm{Pic}}

\title{Fourier-Mukai transforms of curves and principal polarizations}

\author{Marcello Bernardara}

\address{Laboratoire J. A. Dieudonn\'e, Universit\'e
de Nice Sophia Antipolis, Parc Valrose, 06108 Nice Cedex 2 (France).}

\email{bernamar\char`\@math.unice.fr}

\begin{abstract}
Given a Fourier-Mukai transform $\Phi: \derccat \to {\bf{D}}(C')$ between the bounded derived
categories of two smooth projective curves, we verifiy that the
induced map $\phi: J(C) \to J(C')$ between the Jacobian varieties preserves the principal
polarization if and only if $\Phi$ is an equivalence.
\end{abstract}
\maketitle

\section{Introduction}

Let $C$ and $C'$ be two smooth projective curves. Suppose we are given a Fourier-Mukai
functor $\Phi_{\mathcal E}: \derccat \to {\bf{D}}(C')$, with kernel $\mathcal E$, between
the bounded derived categories of the two curves. Recall that $\Phi_{\mathcal E}$ is the
functor
$$\begin{array}{rlll}
\Phi_{\mathcal E} :& \derccat &\longrightarrow &{\bf{D}}(C')\\
& {\mathcal F} &\longmapsto &Rq_{\ast}(p^{\ast}{\mathcal F} \otimes {\mathcal E}),
\end{array}$$
where $p$ and $q$ are the projections of $C \times C'$ respectively on $C$ and $C'$ and
$\mathcal E$ is an object of the bounded derived category ${\bf{D}}(C \times C')$ of the
product of the two curves.\par
If $\Phi_{\mathcal E}$ is an equivalence, it is very well known that we
get an isomorphism between the two curves. We can anyway ask ourselves if this derived
Fourier Mukai transform (DFM for short) does carry an isomorphism between the Jacobian
varieties preserving the principal polarizations.\par
In order to do that, recall the definition of the Jacobian variety as $\pic^0(C)$, the
degree zero part of the Picard group. What we are actually going to do is to make the
DFM descend to an affine map $\Phi_e^P: \pic_{\rational}(C) \to \pic_{\rational}(C')$ between the rational Picard groups.\par
We want to define a morphism $\phi_e^J : J(C) \to J(C')$ compatible with
$\Phi_e^P$. This is given in a unique way, since we know how it acts on the degree zero Picard
group with rational coefficients. All this is done in Section \ref{fromDFM}. In Section
\ref{preservation} we answer the main question and we find a correspondence between the
Torelli Theorem and the derived characterization of a smooth projective curve.
\begin{notation}
For a smooth projective variety $X$ we denote by $\derxcat$ the bounded derived category
of coherent sheaves on $X$.\par
Given an object $\mathcal E$ in $\derxcat$, we denote by $e$ the class $[\mathcal E]$
in $K(X)$ (see section \ref{fromDFM}).\par
We will use the underscore $_{\rational}$ to mean the tensor product with $\rational$.
In particular $\pic_{\rational}(C)$ and $J_{\rational}(C)$ are used to mean respectively
the Picard group and the Jacobian variety with rational coefficients of a smooth projective
curve $C$.\par
Given a product $C \times C'$ of two smooth projective curves, we denote by $p$ and $q$
the projections of $C \times C'$ respectively on $C$ and $C'$.
\end{notation}

\section{From Derived Fourier-Mukai to Jacobian Fourier}\label{fromDFM}

In this section, we describe how a DFM $\Phi_{\mathcal E}$ induces a unique morphism
$\phi_e^J$ between the Jacobian varieties. In order to do that, we describe the morphism
induced by $\Phi_{\mathcal E}$ on the rational Picard group, that is the degree one part
of the rational Chow ring, using standard tecqniques such as Grothendieck-Riemann-Roch
Theorem. What we find is actually an affine map between rational vector
spaces and not a linear morphism. Then the morphism induced on the Jacobian variety
with rational coefficients is the linearization of the affine map restricted to the
degree zero part. Finally, if we know the action of a morphism on the Jacobian with
rational coefficients, there is a unique way to extend it to the Jacobian variety.\par
The way how a DFM between two smooth projective varieties descends at the level of cohomology and gives an isomorphism between the rational cohomology rings is well known. A very good
reference for that is the recent book by Huybrechts \cite{huybrechts}.\par
The first step in making such a descent is going from derived categories to Grothendieck groups. Given $C$ a smooth projective curve, to any object $\mathcal E$ in $\derccat$,
we can associate an element $[{\mathcal E}]$ in the Grothendieck group $K(C)$ by the
alternate sum of the classes of cohomology sheaves of $\mathcal E$.
We thus obtain a map $[\ ]$ from the isomorphism classes of $\derccat$ to the Grothendieck
group $K(C)$.
Given $f: C \to C'$ a projective morphism between smooth projective curves, the pull
back $f^{\ast}: K(C') \to K(C)$ defines a ring homomorphism. The generalized direct image
$f_{!}:K(C) \to K(C')$, defined by $f_{!} {\mathcal F} = \sum (-1)^i R^i f_{\ast}{\mathcal F}$
for any coherent sheaf ${\mathcal F}$ on $C$, defines a group homomorphism.\par
We would like now to define a $K$-theoretic Fourier-Mukai transform (KFM). Given $e$ a class
in $K(C \times C')$, let us define
$$\begin{array}{rlll}
\Phi^K_e :& K(C) &\longrightarrow &K(C')\\
& f &\longmapsto &q_{!}(p^{\ast}f \otimes e).
  \end{array}$$
If we are now given a DFM $\Phi_{\mathcal E} : \derccat \to {\bf{D}}(C')$ with kernel
$\mathcal E$ in ${\bf{D}}(C \times C')$, we obtain the corresponding KFM
$\Phi^K_e : K(C) \to K(C')$ by using the $K$-theoretic kernel $e := [{\mathcal E}]$,
which is a class in  $K(C \times C')$. By the compatibility of $f_!$ and $f^{\ast}$
with $[\ ]$, we easily deduce the following commutative diagram (see \cite{huybrechts},
5.2):
\begin{equation}\label{KandDcompat}
\xymatrix{
\derccat \ar[r]^{\Phi_{\mathcal E}} \ar[d]^{[\ ]} & {\bf{D}}(C') \ar[d]^{[\ ]}\\
K(C) \ar[r]_{\Phi_{[{\mathcal E}]}} & K(C').}
\end{equation}
We want to make a step further and consider rational Chow rings. Consider the exponential Chern character:
$$\begin{array}{rll}
ch:& K(C) \longrightarrow& CH^{\ast}_{\rational}(C)
\end{array}$$
which maps a class of the Grothendieck group to a cycle in the Chow ring with rational coefficients.
For a given $f: C \to C'$ we can define the pull-back $f^{\ast}: CH^{\ast}_{\rational}(C')
 \to CH^{\ast}_{\rational}(C)$ and the direct image $f_{\ast}: CH^{\ast}_{\rational}(C) \to
CH^{\ast}_{\rational}(C')$. Anyway, in order to get a compatibility with the Chern character
$ch$, the Grothendieck-Riemann-Roch Theorem has to be taken into account.\par
\begin{theorem}{\bf{(Grothendieck-Riemann-Roch).}}
Let $f: X \to Y$ a projective morphism of smooth projective varieties. Then for any $e$
in $K(X)$
$$ch(f_{!}(e)) = f_{\ast}(ch(e).Td(f)),$$
where $td(f)$ is the Todd class relative to $f$.
\end{theorem}
In our particular case, the Todd class relative to the
projection $q$ is the Todd class of the relative tangent sheaf:
$$Td(q) = 1 - \frac{1}{2} p^{\ast} K_C,$$
where $K_C$ is the canonical bundle on $C$.\par
Now define the affine map
$$\begin{array}{rlll}
\Phi^P_e :& \pic_{\rational}(C) &\longrightarrow &\pic_{\rational}(C') \\
& M &\longmapsto &q_{\ast}(p^{\ast}M.c_1(e) - \frac{1}{2} c_1(e).p^{\ast}K_C + \frac{1}{2}(c_1^2(e) + 2c_2(e))).
\end{array}$$
The affine map $\Phi^P_e$ is the one induced by the DFM $\Phi_{\mathcal E}$.

\begin{lemma}\label{inductiononpic}
Let $C$ and $C'$ be smooth projective curves and let us suppose there exists a Fourier-Mukai transform
$\Phi_{\mathcal E} : \derccat \to {\bf{D}}(C')$ with kernel $\mathcal E$ in ${\bf{D}}(C \times C')$. Let
$e := [{\mathcal E}]$ be the class in $K(C\times C')$ associated to $\mathcal E$. The
Fourier-Mukai transform $\Phi_{\mathcal E}$ induces the affine map $\Phi_e^P:
\pic_{\rational}(C) \to \pic_{\rational}(C')$ between rational Picard groups,
that means that the diagram
\begin{equation}\label{PandDcompat}
\xymatrix{
\derccat \ar[r]^{\Phi_{\mathcal E}} \ar[d]^{c_1 \circ [\ ]} & {\bf{D}}(C') \ar[d]^{c_1 \circ
[\ ]}\\
\pic_{\rational}(C) \ar[r]_{\Phi_e^P} & \pic_{\rational}(C').}
\end{equation}
is commutative.
\end{lemma}

\begin{proof}
Let us denote by $M$ both an element of the rational Picard group $\pic_{\rational}(C)$
and its class in the Grothendieck group $K(C)$. We want to calculate the first
Chern class $(ch(\Phi_e^K(M))_1$.\par
We have the following chain of equalities:
$$(ch(q_! (p^{\ast} M \otimes e)))_1 = (q_{\ast}(ch(p^{\ast}M \otimes e)(1 - \frac{1}{2}p^{\ast}K_C)))_1,$$
by Grothendieck-Riemann-Roch and $Td(q) = 1 - \frac{1}{2}p^{\ast}K_C$.
$$(q_{\ast}(ch(p^{\ast}M \otimes e).(1 - \frac{1}{2}p^{\ast}K_C)))_1 = q_{\ast}(ch(p^{\ast}M).ch(e).(1 - \frac{1}{2}p^{\ast}K_C))_2.$$
Now let us make it more explicit
\begin{equation}\label{degree2part}
ch(p^{\ast}M).ch(e).(1 - \frac{1}{2}p^{\ast}K_C) = (1 + p^{\ast}M).(r + c_1(e) + \frac{1}{2}(c_1^2(e) + 2c_2(e)).
(1 - \frac{1}{2} p^{\ast}K_C),
\end{equation}
where $r$ is the rank of $e$.\par
We take the degree two part of (\ref{degree2part}), remarking that $p^{\ast}M.p^{\ast}K_C = 0$, and we obtain
\begin{equation}\label{explicitly}
(ch(q_! (p^{\ast} M \otimes e)))_1 = p^{\ast}M.c_1(e) - \frac{1}{2} c_1(e).p^{\ast}K_C + \frac{1}{2}(c_1^2(e) + 2c_2(e)).
\end{equation}
So  the morphism $\Phi_e^P$ between the Picard groups with rational coefficients commutes with the KFM with kernel $e$. Combining this with the commutative diagram (\ref{KandDcompat})
we get the commutative diagram (\ref{PandDcompat}).

\end{proof}
It is clear that the affine map $\Phi_e^P$ restricted to $\pic_{\rational}^0(C)$ does not
give a group morphism to $\pic_{\rational}^0(C')$. Remark anyway that only the first term
of $\Phi_e^P(M)$ depends on $M$, while the other terms are constant with respect
to it.\par
Our aim is to find a morphism between the Jacobian varieties which is compatible with
$\Phi_e^P$ on the degree zero part of the rational Picard group. Such a morphism should
be a group homomorphim first, especially it has to send the zero of $J(C)$ to the zero of $J(C')$.\par
Let us define
$$\begin{array}{rlll}
\phi^J_e :& J(C) &\longrightarrow &J(C') \\
& M &\longmapsto &q_{\ast}(p^{\ast}(M-{{\mathcal O}_C}).c_1(e)),
\end{array}$$
where $e$ is a class in the Grothendieck group $K(C \times C')$ and ${\mathcal O}_C$ is
the unity in $J(C)$. This is the classical Fourier transform with kernel $c_1(e)$ between
the Jacobian varieties, and we are referring to that by JF.\par
Consider the morphism $\phi_e^{J_{\rational}}$, induced by $\phi_e^J$ on
$J_{\rational}(C)$. It coincides with the one induced by $\phi_e^P$
on $J_{\rational}(C)$, which is in turn induced by $\Phi_{\mathcal E}$.
So far we can say that the the DFM with kernel $\mathcal E$ induces on $J_{\rational}(C)$
the morphism $\phi_e^{J_{\rational}}$. This is unique since it is the only
linear morphism induced. We want to make the morphism descend to a morphism
$\phi: J(C) \to J(C')$, but this can be done in a unique way and the result is the JF
$\phi_J^e$.\par
We can then conclude that the JF $\phi_J^e: J(C) \to J(C')$ is the only morphism compatible
with the DFM $\Phi_{\mathcal E}: \derccat \to {\bf{D}}(C')$. Moreover, the correspondence
between DFMs and JFs is functorial.

\begin{lemma}\label{functorial}
The correspondence between derived Fourier-Mukai functors and Fourier transforms on the
Jacobian varieties associating $\phi_e^J$ to $\Phi_{\mathcal E}$ is functorial.
\end{lemma}
\begin{proof}
Given a smooth projective curve $C$, the identity on $\derccat$ is
given by the DFM with kernel ${\mathcal O}_{\Delta}$, the structure sheaf of the diagonal
in $C \times C$. The identity on $J(C)$ clearly corresponds to it.\par
It is clear by Lemma \ref{inductiononpic} that the correspondence between DFMs
and the affine maps is functorial. Indeed if we consider two composable DFMs 
$\Phi_{{\mathcal E}_1}$ and $\Phi_{{\mathcal E}_2}$ and their composition $\Phi_{\mathcal R}$,
the affine maps $\Phi_{e_1}^P$ and $\Phi_{e_2}^P$ are composable and their composition is
given by the affine map $\Phi_r^P$.\par
Now the rational linear map $\phi_e^{J_{\rational}}$ is the linear map associated to
$\Phi_e^P$. Consider in general two composable affine maps $F_1 : V_1 \to V_2$ and
$F_2: V_2 \to V_3$ between vector spaces. Each map $F_i$ induces
the linear map $f_i$. The composition $F_2 \circ F_1$ induces the linear map $f_2 \circ f_1$,
the composition of $f_1$ and $f_2$.\par
This general consideration allows us to state that the correspondence associating the
linear map $\phi_e^{J_{\rational}}$ to the DFM $\Phi_{\mathcal E}$ is functorial.
The last step is just remarking that the functoriality for $\phi_e^{J_{\rational}}$
implies the functoriality for $\phi_e^J$.
\end{proof}

\begin{remark}\label{modifykernel}
Let us observe what happens to the kernel of the JF when we modify the kernel of the DFM.
These remarks will be useful in the next section.\par
The DFM with kernel ${\mathcal E}[1]$ induces the JF with kernel $-c_1(e)$. This can
be computed remarking that $[{\mathcal E}] = -[{\mathcal E}[1]]$ by the definition of
$[\ ]: \derccat \to K(C)$. Here $[1]$ means the shift one step on the right.\par
The DFM with kernel ${\mathcal E}^{\vee}$ induces the JF with kernel $-c_1(e)$. This
can be computed remarking that the Chern polynomial $c_t(e^{\vee})$ satisfies 
$c_t(e^{\vee}) = c_{-t}(e)$.\par
Given line bundles $F$ on $C$ and $F'$ on $C'$, the DFM with kernel ${\mathcal E} \otimes
p^{\ast}F \otimes q^{\ast}(F')$ induces the JF with kernel $c_1(e)$. This can be computed
remarking that $p^{\ast}F.q^{\ast}M = q^{\ast}F'.p^{\ast}M = p^{\ast}F.q^{\ast}F' =0$
for any element $M$ in $J(C)$. In the terminology of \cite{birlange}, Chapter 11, we would
say that we have two equivalent correspondances.
\end{remark}

\section{Preservation of the Principal Polarization}\label{preservation}
So far we know how a DFM $\Phi_{\mathcal E} : \derccat \to {\bf{D}}(C')$ acts on the Jacobians
varieties. The main question is now taken into account: does the morphism
$\phi_e^J$ preserve the principal polarization of $J(C)$?\par

To answer this question, recall that a principal polarization on an abelian variety $A$
defines an isomorphism $\theta_A: A \to \hat{A}$. Given an isogeny $\phi: A \to B$ between two
abelian varieties, we can define the dual isogeny $\hat{\phi}: \hat{B} \to \hat{A}$ between
the dual varieties. If both $A$ and $B$ have a principal polarization, the isogeny $\phi$
respects them if the diagram
$$\xymatrix{
A \ar[r]^{\phi} \ar[d]_{\theta_A} & B \ar[d]^{\theta_B}\\
\hat{A} & \hat{B} \ar[l]^{\hat{\phi}} }$$
is commutative. In the case of the Jacobian variety $J(C)$ of a smooth projective curve $C$ 
we know the principal polarization $\theta_C : J(C) \to
\hat{J}(C)$.\par
We can finally show how to get a positive answer to our question. We identify,
by means of the isomorphisms $\theta_C$ and $\theta_{C'}$ the Jacobian varieties $J(C)$ and
$J(C')$ with their duals. What we have to do is just checking that the composition
$\hat{\phi}_e^J \circ \phi_e^J$ is the identity map on $J(C)$. The dual isomorphism $\hat{\phi}_e^J$ can be obtained as the JF in the
opposite way with the same kernel as $\phi_e^J$. Namely
\begin{equation}\label{conditiononc1e}
\begin{array}{rlll}
\hat{\phi}^J_e : &J(C') &\longrightarrow &J(C)\\
& M' &\longmapsto &p_{\ast} (q^{\ast}(M' - {\mathcal O}_{C'}).c_1(e)),
\end{array}
\end{equation}
see for exemple \cite{birlange}, Chapter 11, Proposition 5.3.\par
Let us come back to derived categories. Given a DFM equivalence
$\Phi_{\mathcal E} : \derccat \to {\bf{D}}(C')$ with kernel $\mathcal E$,
we can describe the kernels ${\mathcal E}_L$ and ${\mathcal E}_R$
of its left and right adjoint. Since $\Phi_{\mathcal E}$ is an equivalence, its adjoints
are its quasi-inverses. The left adjoint of $\Phi_{\mathcal E}$ is the DFM equivlence
$\Phi_{{\mathcal E}_L}: {\bf{D}}(C') \to \derccat$ with kernel
\begin{equation}\label{conditionondfm}
{\mathcal E}_L := {\mathcal E}^{\vee} \otimes q^{\ast} K_{C'} [1].
\end{equation}
We know by remark \ref{modifykernel} that the JF isomorphism induced by the DFM
$\Phi_{{\mathcal E}_L}$ on the Jacobian varieties is given by
$$\begin{array}{rlll}
\phi^J_{e_L} :& J(C') &\longrightarrow &J(C) \\
& M' &\longmapsto &p_{\ast}(q^{\ast}(M'-{{\mathcal O}_{C'}}).c_1(e))
\end{array}$$
Then if $\Phi_{\mathcal E}$ induces on the Jacobian varieties the isomorphism $\phi_e^J$,
its quasi inverse $\Phi_{{\mathcal E}_L}$ induces the dual isomorphism $\hat{\phi}_e^J$.\par 

We can finally state the following Theorem.
\begin{theorem}\label{torelli}
Given two smooth projective curves $C$ and $C'$, a Fourier Mukai functor
$\Phi_{\mathcal E}: \derccat \to {\bf{D}}(C')$ is an equivalence if and only if the
morphism $\phi_e^J: J(C) \to J(C')$ is an isomorphism preserving principal polarization.
\end{theorem}

\begin{proof}
The proof is just a patchwork between everything we have said up to now. Indeed by Remark
\ref{modifykernel}, we can easily check that the DFM with kernel (\ref{conditionondfm})
gives the JF (\ref{conditiononc1e}). The proof follows by Lemma \ref{functorial}.
\end{proof}

\begin{remark}
Recall that for smooth projective curves a derived equivalence always corresponds to an
isomorphism (see for example \cite{huybrechts}, Corollary 5.46).
Theorem \ref{torelli} just states the correspondence between the Torelli Theorem (see for
example \cite{grifhar}, page 359) and the characterization of a curve by its derived
cetegory.\par

This correspondence is no longer valid for higher dimensional smooth projective varieties,
first of all for K3 surfaces. Indeed there exist non isomorphic K3 surfaces which have
equivalent derived categories. This leads to state two distinct Torelli-type Theorems
(see \cite{huytorelli} for an overview). The classical global Torelli Theorem states that
two K3 surfaces are isomorphic if and only if there exists a Hodge isometry between the
degree two cohomology lattices with integer coefficients. The derived global Torelli
Theorem states that two K3 surfaces have equivalent derived categories if and only if
there exist a Hodge isometry between the Mukai lattices.\par
It happens however that, if we deal with smooth projective varieties with ample (anti-)
canonical bundle, an equivalence between derived categories always induces an isomorphism
between the varieties. So we can argue to generalize the correspondence in Theorem
\ref{torelli} to higher dimensional smooth projective varities with ample (anti-) canonical
bundle. A classical example could be the smooth cubic threefold in $\pr^4$, for which we
know the validity of a Torelli-type Theorem.
\end{remark}

\end{document}